\documentclass[aps,prb]{revtex4}

\usepackage{amsmath,amssymb,amsthm}
\usepackage{mathtools}
\usepackage{xspace}
\usepackage[usenames,dvipsnames]{xcolor}
\usepackage[colorlinks,hyperindex,breaklinks]{hyperref}
\usepackage{braket}
\usepackage{hyphenat}
\usepackage{bbold}
\usepackage{enumerate}
\usepackage{tabularx}
\usepackage[hang,flushmargin]{footmisc}

\usepackage[british]{babel}

\usepackage{srcltx}
\usepackage{color}

\usepackage{graphicx}
\usepackage{amsmath}
\usepackage{amssymb}
\usepackage{bm}

\begin{document}

\title{On an infinite number of new families of odd-type Euler sums}
\author{J. Braun}
\affiliation{Ludwig Maximilians-Universit{\"a}t, M{\"u}nchen, Germany}
\author{D. Romberger}
\affiliation{Fakult\"at IV, Abt. BWL, Hochschule Hannover, Germany}
\author{H. J. Bentz}
\affiliation{Institut f\"{u}r Mathematik und Informatik, Universit\"{a}t Hildesheim, Germany}
\date{\today}

\begin{abstract}
We present several sequences of Euler sums involving odd harmonic numbers. The calculational technique is based
on proper two-valued integer functions, which allow to compute these sequences explicitly in terms of zeta values only.
\end{abstract}
\maketitle

\section{Introduction}
Harmonic numbers and their generalizations called hyper-harmonic numbers are defined by
\begin{eqnarray}
H_k = H_k^{(1)} = \sum^{k}_{i=1} \frac{1}{i}~, ~~~~~~~~~~~~  H_k^{(n)} = \sum^{k}_{i=1} \frac{1}{i^n}~.
\end{eqnarray}
It has been discovered by many authors in the past starting with Goldbach 1742 and later on with Euler, Borwein  and
others that linear sums of the type
\begin{eqnarray}
S_{n,m} = \sum^{\infty}_{k=1} \frac{H_k^{(n)}}{k^m} 
\end{eqnarray}
are expressible in terms of zeta values only for an odd weight p=n+m. An excellent overview about these works is found
in the paper of Flajolet and Salvy \cite{fla97} and in citations therein. Further work on similar types of Euler sums can
e found, for example, in \cite {ce17,ce18,ce20}. Although, only few sums of the type
\begin{eqnarray}
T^1_{n,m} = \sum^{\infty}_{k=1} \frac{h_k^{(n)}}{k^m} 
\end{eqnarray}
had been discovered in the past \cite{sit85,zeh07,ade16}. Here the symbol $h_k$ denote odd-type (hyper-)harmonic sums which are defined
as follows:
\begin{eqnarray}
h_k = h_k^{(1)} = \sum^{k}_{i=1} \frac{1}{2i-1}~,~~~~~~~~~~~~  h_k^{(n)} = \sum^{k}_{i=1} \frac{1}{(2i-1)^n}~.
\end{eqnarray}

This way the subject of the present paper is to discover a formalism that allows for an explicit calculation of odd harmonic number
sums in terms of zeta values only for an odd weight p=n+m. In detail we introduce a calculational scheme which is based
on proper two-valued integer functions in correspondence to our last work on Euler sums that was devoted to odd harmonic numbers
and central binomial coefficients \cite{bra20}. We start in the first section with the calculation of a special type of 
nonlinear Euler sums defined as:
\begin{eqnarray}
s(n) = \sum^{\infty}_{k=1} \frac{h_k h^{(2)}_k}{k^{2n+1}}~.
\end{eqnarray}
We show that the corresponding integer function allows for an explicit calculation of sums of the type s(n) in terms of zeta values only.
In the following section we introduce a second type of integer function which in combination with the first type of integer function
allows to compute Euler sums of the type:
\begin{eqnarray}
T^1_{2,m} = \sum^{\infty}_{k=1} \frac{h^{(2)}_k}{k^{2m+1}}~.
\end{eqnarray}
in terms of zeta values only. Furthermore we show that five sub-families $T^2,...,T^6$ which are related to the original sum can also be
computed in terms of zeta values only.   

The last section is devoted to the generalization of our formalism by defining successively higher order integer functions which
allow in a recursive way to compute all Euler sums of the type $T^1_{n,m}$ and the related sub-families for odd weight p=m+n in terms
of zeta values only. This way a complete analogy between the conventional Euler sums $S_{n,m}$ and the odd-type Euler sums $T^1_{n,m}$
can be established.

\section{First family}

Here we start with the first family which is defined as: 
\begin{eqnarray}
s(n) = \sum^{\infty}_{k=1} \frac{h_k h^{(2)}_k}{k^{2n+1}}
\end{eqnarray}
with $n \in \mathbb{N}$. For example for n=1 and n=2 it follows: 
\begin{eqnarray}
s(1) = \sum^{\infty}_{k=1} \frac{h_k h^{(2)}_k}{k^3} = \frac{49}{8} \zeta(3)^2 - \frac{945}{128}\zeta(6)~,
\end{eqnarray}
\begin{eqnarray}
s(2) = \sum^{\infty}_{k=1} \frac{h_k h^{(2)}_k}{k^5} = \frac{651}{8} \zeta(3)\zeta(5) - \frac{343}{16}\zeta(2)\zeta(3)^2 - 
\frac{1575}{32} \zeta(8)
\end{eqnarray}

and so on for all odd powers of k in the denominator. The Euler sum s1(1) has been numerically calculated and compared to
the analytical expression to an accuracy of $10^{-16}$. The numerical value is $s1(1) = 1.3394093155989435$ . For second series it
follows numerically $s2(1)= 1.0567810227967086$.

To be able to compute the corresponding series for all n values we need two identities in form of proper valued integer series.

\subsection{Lemma 1}

The following identities hold:
\begin{eqnarray}
f(k) = \sum^{\infty}_{i=1,i\ne k}\frac{h_i}{i(k-i)} = \sum^{\infty}_{i=1}\frac{h_i}{i(k+i)} - \frac{h_k}{k^2}
-2\frac{h_k^{(2)}}{k}~,
\end{eqnarray}
with
\begin{eqnarray}
\sum^{\infty}_{i=1}\frac{h_i}{i(k+i)} = 2ln(2)\frac{h_k}{k} + \frac{1}{k}
\sum^{k}_{i=1} \frac{H_{i-1}}{2i-1}~.
\end{eqnarray}

\subsection{Proof of Lemma 1}
For k=1 it follows:
\begin{eqnarray}
f(1) = \sum^{\infty}_{i=1,i\ne 1}\frac{h_i}{i(k-i)} = -\sum^{\infty}_{i=2}\frac{h_{i}}{i(i-1)}  = -\sum^{\infty}_{i=1}\frac{h_{i+1}}{i(i+1)}~.
\end{eqnarray}
Analogously we get for k=2:
\begin{eqnarray}
f(2) = \sum^{\infty}_{i=1,i\ne 2}\frac{h_i}{i(k-i)} = -\sum^{1}_{i=1}\frac{h_{i}}{i(i-2)} - \sum^{\infty}_{i=1}\frac{h_{i+2}}{i(i+2)}~.
\end{eqnarray}
Thus it is:
\begin{eqnarray}
f(k) = \sum^{\infty}_{i=1,i\ne k}\frac{h_i}{i(k-i)} = -\sum^{k-1}_{i=1}\frac{h_{i}}{i(i-k)} - \sum^{\infty}_{i=1}\frac{h_{i+k}}{i(i+k)}~.
\end{eqnarray}
Therefore it remains to calculate the two sums on the right side of Eq.~(14). It follows first by partial fraction decomposition:
\begin{eqnarray}
\sum^{\infty}_{i=1}\frac{h_{i+k}}{i(i+k)} = \frac{1}{k}\sum^{\infty}_{i=1}\frac{h_{i+k}}{i} -\frac{1}{k} \sum^{\infty}_{i=1}\frac{h_{i}}{i}
+\sum^{k}_{i=1}\frac{h_{i}}{i}~.
\end{eqnarray}
With
\begin{eqnarray}
\sum^{\infty}_{i=1}\frac{h_{i+k}}{i} = \sum^{\infty}_{i=1}\frac{h_i}{i} + \sum^{k}_{m=1} \sum^{\infty}_{i=1}\frac{1}{i(2m+2i-1)} =
\sum^{\infty}_{i=1}\frac{h_{i}}{i} + \sum^{k}_{m=1} \left(\frac{2h_m}{2m-1} - \frac{2ln(2)}{2m-1} \right)
\end{eqnarray}
it results
\begin{eqnarray}
\sum^{\infty}_{i=1}\frac{h_{i+k}}{i(i+k)} = -2ln(2)\frac{h_k}{k} + \frac{1}{k} \sum^{k}_{i=1}\frac{h_{i}}{i} + \frac{2}{k}
\sum^{k}_{i=1}\frac{h_i}{2i-1}~,
\end{eqnarray}
or
\begin{eqnarray}
\sum^{\infty}_{i=1}\frac{h_{i+k}}{i(i+k)} = -2ln(2)\frac{h_k}{k} + \frac{1}{k} \sum^{k}_{i=1}\frac{h_i}{i} + \frac{h_k^{(2)}}{k}
+ \frac{h_k^2}{k}~.
\end{eqnarray}
For the second sum it follows first:
\begin{eqnarray}
\sum^{k-1}_{i=1}\frac{h_{i}}{i(k-i)} = \frac{1}{k}\sum^{k-1}_{i=1}\frac{h_i}{i} + \frac{1}{k}\sum^{k-1}_{i=1}\frac{h_i}{k-i}~. 
\end{eqnarray}
With
\begin{eqnarray}
\sum^{k-1}_{i=1}\frac{h_i}{k-i} = \sum^{i-1}_{k=1}\frac{h_{i-k}}{k} = F(i)
\end{eqnarray}
we get
\begin{eqnarray}
F(i+1) = \sum^{i-1}_{k=1}\frac{h_{i-k}}{k} + \sum^{i}_{k=1} \frac{1}{k(2i-2k+1)} = F(i) + \sum^{i}_{k=1} \frac{1}{k(2i-2k+1)}~.
\end{eqnarray}
Thus an inhomogeneous difference equation of first oder results 
\begin{eqnarray}
F(i+1) - F(i) = \frac{1}{2i+1} \left( H_i + 2h_i \right)
\end{eqnarray}
with the solution
\begin{eqnarray}
F(i) = \sum^{i}_{k=1} \frac{H_{k-1}}{2k-1} +2\sum^{i}_{k=1}\frac{h_{k}}{2k-1} -2\sum^{i}_{k=1}\frac{1}{(2k-1)^2}~.  
\end{eqnarray}
Using the identity (3.94) calculated by Adegoke \cite{ade16} it follows:
\begin{eqnarray}
\sum^{k-1}_{i=1}\frac{h_i}{k-i} = H_k h_k - \sum^{k}_{i=1} \frac{h_i}{i} + h_k^{(2)} + h_k^2~.
\end{eqnarray}
Therefore it follows for the second sum:
\begin{eqnarray}
\sum^{k-1}_{i=1}\frac{h_{i}}{i(k-i)} = \frac{H_k h_k}{k} -\frac{h_k}{k^2} + \frac{h_k^{(2)}}{k} + \frac{h_k^2}{k}~. 
\end{eqnarray}
This way Eq.~(14) results in:
\begin{eqnarray}
\sum^{\infty}_{i=1,i\ne k}\frac{h_i}{i(k-i)} = 2ln(2)\frac{h_k}{k} - \frac{h_k}{k^2} - 2\frac{h_k^{(2)}}{k} + \frac{H_k h_k}{k}
-\frac{1}{k} \sum^{k}_{i=1} \frac{h_i}{i}~.
\end{eqnarray} 
It remains to proof Eq.~(11). It follows first from partial fraction decomposition
\begin{eqnarray}
\sum^{\infty}_{i=1}\frac{h_i}{i(i+k)} = \frac{1}{k}\sum^{k}_{i=1}\frac{h_i}{i} + 2ln(2) \frac{h_k}{k} + \frac{H_k h_k}{k}
-\frac{2}{k}\sum^{k}_{i=1}\frac{h_i}{2i+2k-1}~.
\end{eqnarray}
With
\begin{eqnarray}
\sum^{k}_{i=1}\frac{h_i}{2i+2k-1} = \sum^{k}_{i=1}\frac{h_{k+1-i}}{2i-1} = F(k) 
\end{eqnarray}
it follows
\begin{eqnarray}
F(k) - F(k-1) = \sum^{k}_{i=1}\frac{1}{(2i-1)(2i+2k-1)} = \frac{h_k}{k}~. 
\end{eqnarray}
The solution of the corresponding inhomogeneous difference equation is
\begin{eqnarray}
F(k) = \sum^{k}_{i=1} \frac{h_i}{i}~.
\end{eqnarray}
Therefore, we get  
\begin{eqnarray}
\sum^{\infty}_{i=1}\frac{h_i}{i(i+k)} = 2ln(2) \frac{h_k}{k} + \frac{H_k h_k}{k} -\frac{1}{k}\sum^{k}_{i=1}\frac{h_i}{i} ~,
\end{eqnarray}
and with the identity (3.10) calculated by Adegoke \cite{ade16} Eq.~(11) results. Thus lemma 1 is proofed.

The first family of nonlinear Euler sums follows now by calculating the following expression for $n \in \mathbb{N}$
\begin{eqnarray}
s(n) = \sum^{\infty}_{k=1}\frac{h_k}{k^{2n}}f(k)~. 
\end{eqnarray}
As an example we compute for n=1
\begin{eqnarray}
\sum^{\infty}_{k=1}\frac{h_k}{k^{2}}f(k) = \sum^{\infty}_{k=1} \frac{h_k}{k^2} \sum^{\infty}_{i=1}\frac{h_i}{i(k+i)} -
\sum^{\infty}_{k=1} \frac{h_k^2}{k^4} - 2\sum^{\infty}_{k=1} \frac{h_k h_k^{(2)}}{k^3} .
\end{eqnarray}
With
\begin{eqnarray}
\sum^{\infty}_{k=1} \frac{h_k}{k^2} \sum^{\infty}_{i=1}\frac{h_i}{i(k+i)} = \sum^{\infty}_{i=1} \frac{h_i}{i}
\sum^{\infty}_{k=1}\frac{h_k}{k^2(k+i)} \nonumber
\end{eqnarray}
\begin{eqnarray}
= \sum^{\infty}_{i=1} \frac{h_i}{i^2} \sum^{\infty}_{k=1}\frac{h_k}{k^2} - \sum^{\infty}_{i=1} \frac{h_i}{i^2} 
\sum^{\infty}_{k=1}\frac{h_k}{k(i+k)}~.
\end{eqnarray}
Thus it follows:
\begin{eqnarray}
\sum^{\infty}_{k=1} \frac{h_k}{k^2} \sum^{\infty}_{i=1}\frac{h_i}{i(k+i)} = 
\frac{1}{2} \left(\sum^{\infty}_{k=1}\frac{h_k}{k^2} \right)^2~.
\end{eqnarray}
As a next step we compute
\begin{eqnarray}
\sum^{\infty}_{k=1} \frac{h_k}{k^2}\sum^{\infty}_{i=1,i\ne k}\frac{h_i}{i(k-i)} =  -\sum^{\infty}_{i=1}
\frac{h_i}{i^2} \sum^{\infty}_{k=1,k\ne i}\frac{h_k}{k^2} -
\sum^{\infty}_{i=1} \frac{h_i}{i^2} \sum^{\infty}_{k=1,k\ne i}\frac{h_k}{k(i-k)}~.
\end{eqnarray}
Thus it follows:
\begin{eqnarray}
\sum^{\infty}_{k=1} \frac{h_k}{k^2}\sum^{\infty}_{i=1,i\ne k}\frac{h_i}{i(k-i)} = - \left(\sum^{\infty}_{k=1}\frac{h_k}{k^2} \right)^2 +
\sum^{\infty}_{k=1}\frac{h_k}{k^4}~.
\end{eqnarray}
Inserting Eq.~(35) and Eq.~(37) in Eq.~(33) we get:
\begin{eqnarray}
\sum^{\infty}_{k=1} \frac{h_k h_k^{(2)}}{k^3} = \frac{1}{2} \left(\sum^{\infty}_{k=1}\frac{h_k}{k^2} \right)^2 - \frac{3}{2}
\sum^{\infty}_{k=1}\frac{h_k}{k^4}~.
\end{eqnarray}
Here, both Euler sums on the right side are known explicitly in terms of zeta values, and this way the sum of both terms
results in Eq.~(8).

\section{Second family type 1}

The second family of type 1 is defined as follows: 
\begin{eqnarray}
T^1_{2,m} = \sum^{\infty}_{k=1} \frac{h^{(2)}_k}{k^{2m+1}}
\end{eqnarray}
with $n \in \mathbb{N}$. For example for m=1 and m=2 it follows: 
\begin{eqnarray}
\sum^{\infty}_{k=1} \frac{h^{(2)}_k}{k^3} = \frac{35}{4} \zeta(2)\zeta(3) - \frac{31}{2} \zeta(5)
\end{eqnarray}
and
\begin{eqnarray}
\sum^{\infty}_{k=1} \frac{h^{(2)}_k}{k^5} = \frac{217}{4} \zeta(2)\zeta(5) - \frac{7}{2}\zeta(3)\zeta(4) - \frac{381}{4} \zeta(7)~,
\end{eqnarray}
and so on for all odd powers of k in the denominator. The first series is known from literature \cite{zeh07}. For the second sum it
follows numerically $T^1_{2,2} = 1.01413007995319209 $ with an accuracy of $10^{-16}$ compared to the analytical expression.
To be able to compute the corresponding series for all n values we again need a new identity in form of a proper valued integer series.

\subsection{Lemma 2}

The following identity holds:
\begin{eqnarray}
g_{2n}(k) = \sum^{\infty}_{i=1,i\ne k}\frac{1}{i^{2n}(k-i)} = \frac{H_k}{k^{2n}} - \frac{2n+1}{k^{2n+1}} 
+ \sum^{2n-1}_{i=1}\zeta(2n+1-i)\frac{1}{k^i}~. 
\end{eqnarray}

\subsection{Proof of Lemma 2}
For k=1 it follows:
\begin{eqnarray}
g_{2n}(1) = \sum^{\infty}_{i=1,i\ne 1}\frac{1}{i^{2n}(1-i)} = -\sum^{\infty}_{i=2}\frac{1}{i^{2n}(i-1)}
= -\sum^{\infty}_{i=1}\frac{1}{i(i+1)^{2n}}~.
\end{eqnarray}
Analogously we get for k=2:
\begin{eqnarray}
g_{2n}(2) = \sum^{\infty}_{i=1,i\ne 2}\frac{1}{i^{2n}(2-i)} = -\sum^{1}_{i=1}\frac{1}{i^{2n}(i-2)}
- \sum^{\infty}_{i=1}\frac{1}{i(i+2)^{2n}}~.
\end{eqnarray}
Thus it is:
\begin{eqnarray}
g_{2n}(k) = \sum^{\infty}_{i=1, i\ne n}\frac{1}{i^{2n}(k-i)} = -\sum^{k-1}_{i=1}\frac{1}{i^{2n}(i-k)} 
- \sum^{\infty}_{i=1}\frac{1}{i(i+k)^{2n}}~.
\end{eqnarray}
For the second sum on the right side of Eq.~(45) we simply get by partial fraction decomposition:
\begin{eqnarray}
\sum^{\infty}_{i=1}\frac{1}{i(i+k)^{2n}} = \frac{H_k}{k^{2n}} - \sum^{2n-1}_{i=1}\zeta(2n+1-i)\frac{1}{k^i}
+ \sum^{2n-1}_{i=1}H_k^{(2n+1-i)}\frac{1}{k^i}~.
\end{eqnarray}
It remains to calculate the finite sum:
\begin{eqnarray}
\sum^{k-1}_{i=1}\frac{1}{i^{2n}(k-i)} = \sum^{2n}_{i=1}\frac{H_{k-1}^{(2n+1-i)}}{k^i} + \frac{1}{k^{2n}}
\sum^{k-1}_{i=1}\frac{1}{k-i}~.
\end{eqnarray}
This result again follows from partial fraction decomposition. With
\begin{eqnarray}
\sum^{k-1}_{i=1}\frac{1}{k-i} = H_{k-1}
\end{eqnarray}
we get
\begin{eqnarray}
\sum^{k-1}_{i=1}\frac{1}{i^{2n}(k-i)} = \sum^{2n}_{i=1}\frac{H_{k-1}^{2n+1-i}}{k^i} + \frac{H_{k-1}}{k^{2n}}~.
\end{eqnarray}
Addition of both sums gives Eq.~(42). Thus the lemma is proofed.

When we now calculate the following expression for $n \in \mathbb{N}$
\begin{eqnarray}
T^1_{2,m} = \sum^{\infty}_{k=1}\frac{h_k}{k}g_{2m}(k)
\end{eqnarray}
the second family of nonlinear Euler sums follows. As an example we compute the series for m=1. We get: 
\begin{eqnarray}
\sum^{\infty}_{k=1}\frac{h_k}{k}g_2(k) = \sum^{\infty}_{k=1} \frac{h_k}{k} \left(\frac{H_k}{k^2}-
\frac{3}{k^3}+\zeta(2)\frac{1}{k} \right)  = \sum^{\infty}_{k=1} \frac{h_k h^{(2)}_k}{k^3} -3\sum^{\infty}_{k=1} \frac{h_k}{k^4}
+\zeta(2)\frac{h_k}{k^2}~, 
\end{eqnarray}
and 
\begin{eqnarray}
\sum^{\infty}_{k=1}\frac{h_k}{k}g_2(k) = \sum^{\infty}_{i=1}\frac{1}{i^2} \sum^{\infty}_{k=1,k\ne i}\frac{h_k}{k(k-i)} =
-\sum^{\infty}_{k=1}\frac{h_k}{k^2}\sum^{\infty}_{i=1}\frac{1}{i^2(k+i)} + \sum^{\infty}_{k=1}\frac{h_k}{k^4}
+2\sum^{\infty}_{k=1}\frac{h_k^{(2)}}{k^3}~.
\end{eqnarray}
Thus it follows:
\begin{eqnarray}
\sum^{\infty}_{k=1}\frac{h_k^{(2)}}{k^3} = \zeta(2)\sum^{\infty}_{k=1}\frac{h_k}{k^2} -2 \sum^{\infty}_{k=1}\frac{h_k}{k^4}~.
\end{eqnarray}
Here, both Euler sums on the right side are known explicitly in terms of zeta values, and this way the sum of both terms
results in Eq.~(40) .
Analogously Eq.~(41) follows from
\begin{eqnarray}
\sum^{\infty}_{k=1}\frac{h_k^{(2)}}{k^5} = \zeta(4)\sum^{\infty}_{k=1}\frac{h_k}{k^2} + \zeta(2)\sum^{\infty}_{k=1}\frac{h_k}{k^4} -
3\sum^{\infty}_{k=1}\frac{h_k}{k^6}~,
\end{eqnarray}
where again all Euler sums on the right side are known explicitly in terms of zeta values.

\section{Second family type 2}

The second family of type 2 is defined as follows: 
\begin{eqnarray}
T^2_{2,m} = \sum^{\infty}_{k=1} \frac{H^{(2m+1)}_k}{(2k-1)^2}~.
\end{eqnarray}
For m=1 and m=2 it follows: 
\begin{eqnarray}
\sum^{\infty}_{k=1} \frac{H^{(3)}_k}{(2k-1)^2} =  \frac{31}{2}\zeta(5) - 8 \zeta(2)\zeta(3) +\zeta(3) + 10\zeta(2) -24 ln(2)~, 
\end{eqnarray}
and
\begin{eqnarray}
\sum^{\infty}_{k=1} \frac{H^{(5)}_k}{(2k-1)^2} =  \frac{381}{4}\zeta(7) - \frac{107}{2} \zeta(2)\zeta(5) - \frac{7}{2}\zeta(3)\zeta(4) 
+\zeta(5) + 4\zeta(4) + 12 \zeta(3) + 56\zeta(2) -160 ln(2)~.                           
\end{eqnarray}
and so on for all odd powers of k in the nominator. These expressions follow from the well known identity (3.94) calculated by
Adegoke \cite{ade16} . 

\section{Second family type 3}

The second family of type 3 is defined as follows: 
\begin{eqnarray}
T^3_{2,m} = \sum^{\infty}_{k=1} \frac{h^{(2)}_k}{(2k-1)^{2m+1}}~.
\end{eqnarray}
For m=1 and m=2 it follows: 
\begin{eqnarray}
\sum^{\infty}_{k=1} \frac{h^{(2)}_k}{(2k-1)^3} = \frac{31}{64}\zeta(5) + \frac{9}{32}\zeta(2)\zeta(3)~,
\end{eqnarray}
and
\begin{eqnarray}
\sum^{\infty}_{k=1} \frac{h^{(2)}_k}{(2k-1)^5} = \frac{127}{256}\zeta(7) + \frac{15}{128}\zeta(2)\zeta(5) + \frac{15}{64}\zeta(3)\zeta(4)~.
\end{eqnarray}

This family of nonlinear Euler sums follows by calculating the function r$(n)$ for n$\in \mathbb{N}$: 
\begin{eqnarray}
r(n) = \sum^{\infty}_{i=1} \frac{1}{i} \sum^{\infty}_{k=1} \frac{1}{k^{2n}(i+2k)^2}~. 
\end{eqnarray}

As an example we compute Eq.~(59) for n=1. It follows first by use of partial fraction decomposition:
\begin{eqnarray}
\sum^{\infty}_{i=1} \frac{1}{i} \sum^{\infty}_{k=1} \frac{1}{k^{2}(i+2k)^2} =
\sum^{\infty}_{i=1} \frac{1}{i^3} \sum^{\infty}_{k=1} \frac{1}{k^2} -4\sum^{\infty}_{i=1} \frac{1}{i^3} 
\sum^{\infty}_{k=1} \frac{1}{k(i+2k)} +4\sum^{\infty}_{i=1} \frac{1}{i^3} \sum^{\infty}_{k=1} \frac{1}{(i+2k)^2}~.
\end{eqnarray}
Splitting of the sums in odd and even contributions we get:
\begin{eqnarray}
\sum^{\infty}_{i=1} \frac{1}{i} \sum^{\infty}_{k=1} \frac{1}{k^{2}(i+2k)^2} &=&
\zeta(2)\zeta(3) - 4\sum^{\infty}_{i=1}\frac{1}{(2i-1)^3}\sum^{\infty}_{k=1}\frac{1}{k(2i+2k-1)} -
\frac{1}{4}\sum^{\infty}_{i=1}\frac{1}{i^3}\sum^{\infty}_{k=1}\frac{1}{k(i+k)} \nonumber \\ &+&
\sum^{\infty}_{i=1}\frac{1}{(2i-1)^3}\sum^{\infty}_{k=1}\frac{1}{(2i+2k-1)^2} +
\frac{1}{8} \sum^{\infty}_{i=1}\frac{1}{i^3}\sum^{\infty}_{k=1}\frac{1}{(i+k)^2}~,
\end{eqnarray}
where we have used the identities
\begin{eqnarray}
\sum^{\infty}_{k=1}\frac{1}{k(2i+2k-1)} = \frac{2}{2i-1} \left(h_i - ln(2) \right)~,
\end{eqnarray}
and
\begin{eqnarray}
\sum^{\infty}_{k=1}\frac{1}{k(i+k)} = \frac{H_i}{i}~.
\end{eqnarray}
Eq.~(64) and Eq.~(65) can be simply proofed by partial fraction decomposition. It follows then:
\begin{eqnarray}
\sum^{\infty}_{i=1} \frac{1}{i} \sum^{\infty}_{k=1} \frac{1}{k^{2}(i+2k)^2} &=&
\frac{9}{8}\zeta(2)\zeta(3)-8\sum^{\infty}_{i=1} \frac{h_i}{(2i-1)^4} + 8ln(2) \sum^{\infty}_{i=1} \frac{1}{(2i-1)^4}-
\sum^{\infty}_{i=1}\frac{H_i}{i^4} + 3\zeta(2)\sum^{\infty}_{i=1} \frac{1}{(2i-1)^3} \nonumber \\ &-& 
\frac{1}{8} \sum^{\infty}_{i=1}\frac{H_i^{(2)}}{i^3} - 4\sum^{\infty}_{i=1} \frac{h_i^{(2)}}{(2i-1)^3}~.
\end{eqnarray}

Finally we get:
\begin{eqnarray}
\sum^{\infty}_{i=1} \frac{1}{i} \sum^{\infty}_{k=1} \frac{1}{k^{2}(i+2k)^2} = -\frac{65}{16}\zeta(5)+\frac{35}{8}\zeta(2)\zeta(3)
-4 \sum^{\infty}_{i=1} \frac{h_i^{(2)}}{(2i-1)^3}~.
\end{eqnarray}

Furthermore, it holds:
\begin{eqnarray}
\sum^{\infty}_{i=1} \frac{1}{i(i+2k)^2} = \frac{h_k}{4k^2}+\frac{H_k}{8k^2}-\zeta(2)\frac{1}{2k}+\frac{h_k^{(2)}}{2k}
+\frac{H_k^{(2)}}{8k}~.
\end{eqnarray}
Again Eq.~(68) can be simply proofed by partial fraction decomposition. Thus it follows:
\begin{eqnarray}
\sum^{\infty}_{k=1} \frac{1}{k^2} \sum^{\infty}_{i=1} \frac{1}{i(i+2k)^2} = \frac{1}{4} \sum^{\infty}_{k=1}\frac{h_k}{k^4} +
\frac{1}{8} \sum^{\infty}_{k=1}\frac{H_k}{k^4} - \frac{1}{2}\zeta(2)\sum^{\infty}_{k=1}\frac{1}{k^3} +
\frac{1}{2} \sum^{\infty}_{k=1}\frac{h_k^{(2)}}{k^3} + \frac{1}{8} \sum^{\infty}_{k=1}\frac{H_k^{(2)}}{k^3} .
\end{eqnarray}
Here, all Euler sums on the right side are known explicitly in terms of zeta values, and this way the total sum results in: 
\begin{eqnarray}
\sum^{\infty}_{k=1} \frac{1}{k^2} \sum^{\infty}_{i=1} \frac{1}{i(i+2k)^2} = \frac{26}{8}\zeta(2)\zeta(3)-6\zeta(5)~.
\end{eqnarray}
As the left sides of Eq.~(67) and Eq.~(69) are equal we get Eq.~(59) .

\section{Second family typs 4}

The second family of type 4 is defined as follows: 
\begin{eqnarray}
T^4_{2,m} = \sum^{\infty}_{k=1} \frac{h^{(2m+1)}_k}{(2k-1)^2}~.
\end{eqnarray}
For m=1 and m=2 the result is: 
\begin{eqnarray}
\sum^{\infty}_{k=1} \frac{h^{(3)}_k}{(2k-1)^2} = \frac{31}{64}\zeta(5) + \frac{3}{8}\zeta(2)\zeta(3)~,
\end{eqnarray}
and
\begin{eqnarray}
\sum^{\infty}_{k=1} \frac{h^{(5)}_k}{(2k-1)^2} = \frac{127}{256}\zeta(7) + \frac{39}{64}\zeta(2)\zeta(5) - \frac{15}{64}\zeta(3)\zeta(4)
\end{eqnarray}
and so on for all odd powers of k in the nominator. These expressions follow from the well known identity (3.85) given by Adegoke \cite{ade16}. 

\section{Second family type 5}

The second family of type 5 is defined as follows: 
\begin{eqnarray}
T^5_{2,m} = T^1_{n,2} = \sum^{\infty}_{k=1} \frac{h^{(2m+1)}_k}{k^2}~.
\end{eqnarray}
For m=1 and m=2 the result is: 
\begin{eqnarray}
\sum^{\infty}_{k=1} \frac{h^{(3)}_k}{k^2} = \frac{93}{8} \zeta(5) - \frac{21}{4} \zeta(2)\zeta(3)~,
\end{eqnarray}
and
\begin{eqnarray}
\sum^{\infty}_{k=1} \frac{h^{(5)}_k}{k^2} = \frac{1905}{64} \zeta(7) - \frac{93}{8} \zeta(2)\zeta(5) - \frac{105}{16}\zeta(3)\zeta(4)
\end{eqnarray}
and so on for all odd powers of k in the nominator.

The series belonging to the fifth family of nonlinear Euler sums can be calculated by use of the following expression:
\begin{eqnarray}
\sum^{\infty}_{k=1} \frac{H_{k}^{(2n+1)}}{k^2} = \sum^{\infty}_{k=1} \frac{H_{2k-1}^{(2n+1)}}{(2k-1)^2} +
\frac{1}{4}\frac{H_{2k}^{(2n+1)}}{k^2}~.
\end{eqnarray} 
As an example we compute Eq.~(75). For n=1 we get: 
\begin{eqnarray}
\sum^{\infty}_{k=1} \frac{H_{k}^{(3)}}{k^2} = \sum^{\infty}_{k=1} \frac{H_{2k-1}^{(3)}}{(2k-1)^2} +
\frac{1}{4}\frac{H_{2k}^{(3)}}{k^2} = \frac{11}{2}\zeta(5) - 2\zeta(2)\zeta(3)~. 
\end{eqnarray}
The Euler sum on the left side is known from literature \cite{zeh07}. With 
\begin{eqnarray}
H_{2k}^{(3)} = h_k^{(3)} + \frac{1}{8} H_{k}^{(3)}
\end{eqnarray} 
it results:
\begin{eqnarray}
\frac{11}{2}\zeta(5) - 2\zeta(2)\zeta(3) = \sum^{\infty}_{k=1} \frac{H_{2k}^{(3)}}{(2k-1)^2} + \frac{1}{4} \sum^{\infty}_{k=1} 
\frac{H_{2k}^{(3)}}{k^2} - \frac{1}{8} \sum^{\infty}_{k=1} \frac{1}{k^3(2k-1)^2} \nonumber
\end{eqnarray} 
\begin{eqnarray}
= \sum^{\infty}_{k=1} \frac{h_{k}^{(3)}}{(2k-1)^2} + \frac{1}{8} \sum^{\infty}_{k=1} \frac{H_{k}^{(3)}}{(2k-1)^2} +
\sum^{\infty}_{k=1} \frac{h_{k}^{(3)}}{k^2} + \frac{1}{32} \sum^{\infty}_{k=1} \frac{H_{k}^{(3)}}{k^2} -
\frac{1}{8} \sum^{\infty}_{k=1} \frac{1}{k^3(2k-1)^2}~.
\end{eqnarray} 
Thus it follows:
\begin{eqnarray}
\frac{h_{k}^{(3)}}{k^2} = \frac{11}{2}\zeta(5) - 2\zeta(2)\zeta(3) - \sum^{\infty}_{k=1} \frac{h_{k}^{(3)}}{(2k-1)^2} -
\frac{1}{8} \sum^{\infty}_{k=1} \frac{H_{k}^{(3)}}{(2k-1)^2} - \frac{1}{32} \sum^{\infty}_{k=1} \frac{H_{k}^{(3)}}{k^2} +
\frac{1}{8} \sum^{\infty}_{k=1} \frac{1}{k^3(2k-1)^2}~.
\end{eqnarray} 
All sums on the right side are known. This way we get Eq.~(75).

\section{Second family type 6}

The second family of type 6 is defined as follows: 
\begin{eqnarray}
T^6_{2,m} = \sum^{\infty}_{k=1} \frac{H^{(2)}_k}{(2k-1)^{2m+1}}~.
\end{eqnarray}
For m=1 and m=2 we get: 
\begin{eqnarray}
\sum^{\infty}_{k=1} \frac{H^{(2)}_k}{(2k-1)^3} = \frac{49}{8} \zeta(2)\zeta(3) - \frac{93}{8}\zeta(5)
+\frac{7}{2}\zeta(3) - 7\zeta(2) + 12 ln(2)~,                           
\end{eqnarray}
and
\begin{eqnarray}
\sum^{\infty}_{k=1} \frac{H^{(2)}_k}{(2k-1)^5} = \frac{403}{32} \zeta(2)\zeta(5) + \frac{105}{16}\zeta(3)\zeta(4) - \frac{1905}{64}\zeta(7)
+ \frac{31}{8} \zeta(5) -\frac{15}{2}\zeta(4) +\frac{21}{2}\zeta(3) -13\zeta(2) +20 ln(2)                           
\end{eqnarray}
and so on for all odd powers of k in the denominator. These expressions again follow from the well known identity (3.94) given by
Adegoke \cite{ade16} .

\section{Generalization to higher n and m values}

The calculational scheme introduced above can be used to compute recursively all higher order Euler sums of the type: 
\begin{eqnarray}
T^1_{n,m} = \sum^{\infty}_{k=1} \frac{h_k^{(n)}}{k^m}~,
\end{eqnarray}
as well as all sub-families $T^2_{n,m}, ... , T^6_{n,m}$ where p=n+m always is of odd weight with $m,n \in \mathbb{N}$. 

The formalism starts with choosing n=3 for an arbitrarily given m value in $T^1_{n,m}$. For example for n=3 and m=4 it follows: 
\begin{eqnarray}
T^1_{3,4} = \sum^{\infty}_{k=1} \frac{h_k^{(3)}}{k^4} = \frac{1905}{16} \zeta(7) - \frac{279}{4} \zeta(2)\zeta(5)~.
\end{eqnarray}
Furthermore, for n=3 and m=6 we get:
\begin{eqnarray}
T^1_{3,6} = \sum^{\infty}_{k=1} \frac{h_k^{(3)}}{k^6} = \frac{3}{4}\zeta(2) \sum^{\infty}_{k=1}\frac{h_k}{k^6} + 
\frac{1}{2}\zeta(4) \sum^{\infty}_{k=1}\frac{h_k^{(2)}}{k^3} +
\frac{1}{2}\zeta(2) \sum^{\infty}_{k=1}\frac{h_k^{(2)}}{k^5} + 
\frac{7}{2} \sum^{\infty}_{k=1}\frac{h_k^{(2)}}{k^7}~.
\end{eqnarray}
All sums on the right side are known, where sums of type $T^1_{2,m}$ appear. This fact indicates that the calculation of higher order
T's has to be done recursively but nevertheless is possible in terms of zeta values only. It follows numerically
$T^1_{3,4} = 1.08556003490415209$ with an accuracy of $10^{-16}$ compared to the analytical expression and $T^1_{3,6} = 1.01800033232122239$.
To be able to compute the corresponding series for all n values we again need a new identity in form of a proper valued integer series.
The procedure can be applied for all
even powers of k in the denominator in an analogous way. For an explicit calculation a new two-valued integer function is needed.

\subsection{Lemma 3}

The following identities hold:
\begin{eqnarray}
f_{2n-1}(k) = \sum^{\infty}_{i=1,i\ne k}\frac{h^{(2n-1)}_i}{i(k-i)} = \sum^{\infty}_{i=1}\frac{h^{(2n-1)}_i}{i(k+i)}
- (4n-2) \frac{h^{(2n)}_k}{k} - \frac{h^{(2n-1)}_k}{k^2} + 4 \sum_{i=1}^{n-1} \left( 1-\frac{1}{2^{2i}} \right)\zeta(2i)
\frac{h^{(2n-2i)}_k}{k}
\end{eqnarray}
and
\begin{eqnarray}
f_{2n}(k) = \sum^{\infty}_{i=1,i\ne k}\frac{h^{(2n)}_i}{i(k-i)} = -\sum^{\infty}_{i=1}\frac{h^{(2n)}_i}{i(k+i)}
- 4n\frac{h^{(2n+1)}_k}{k} - \frac{h^{(2n)}_k}{k^2} + 4 \sum_{i=1}^n\left(1-\frac{1}{2^{2i}} \right)\zeta(2i)
\frac{h^{(2n-2i+1)}_k}{k}~.
\end{eqnarray}
Compared to the identity defined in lemma 1 we see that instead of $h_i$ the hyper-harmonic function $h^{(n)}_i$ appears in the
corresponding identities.

\subsection{Proof of Lemma 3}
In analogy to Eq.~(10) it follows first:
\begin{eqnarray}
f_n(k) = \sum^{\infty}_{i=1,i\ne k}\frac{h^{(n)}_i}{i(k-i)} = -\sum^{k-1}_{i=1}\frac{h^{(n)}_{i}}{i(i-k)} - 
\sum^{\infty}_{i=1}\frac{h^{(n)}_{i+k}}{i(i+k)}~.
\end{eqnarray}
For the infinite sum on the right side we get:
\begin{eqnarray}
\sum^{\infty}_{i=1}\frac{h^{(n)}_{i+k}}{i(i+k)} = \frac{1}{k}\sum^{\infty}_{i=1}\frac{h^{(n)}_{i+k}}{i} -
\frac{1}{k}\sum^{\infty}_{i=1}\frac{h^{(n)}_{i+k}}{i+k} = \frac{1}{k} \sum^{k}_{i=1}\frac{h^{(n)}_i}{i} +
\frac{1}{k} \sum^{k}_{m=1} \left( \sum^{\infty}_{i=1} \frac{1}{i (2i+2m-1)^n} \right)~.
\end{eqnarray}
As a next step we calculate the finite sum
\begin{eqnarray}
\sum^{k-1}_{i=1}\frac{h^{(n)}_{i}}{i(i-k)} = \frac{1}{k}\sum^{k-1}_{i=1}\frac{h^{(n)}_i}{i} + \frac{1}{k} 
\sum^{k-1}_{i=1}\frac{h^{(n)}_i}{k-i}~.
\end{eqnarray}
With
\begin{eqnarray}
\sum^{k-1}_{i=1}\frac{h^{(n)}_i}{k-i} = \sum^{i-1}_{k=1}\frac{h^{(n)}_{i-k}}{k} = F(i)
\end{eqnarray}
it follows:
\begin{eqnarray}
F(m+1) - F(m) = \sum^{m}_{i=1} \frac{1}{i(2i-2m+1)^n}~.
\end{eqnarray}
Thus the solution is
\begin{eqnarray}
F(k+1) = \sum^{k}_{m=1} \left( \sum^{m}_{i=1} \frac{1}{i(2i-2m+1)^n} \right)~.
\end{eqnarray}
By partial fraction decomposition of the inner sum it results:
\begin{eqnarray}
F(k+1) = \sum^{k}_{m=1} \frac{H_m}{(2m+1)^n} +2\sum^{k}_{m=1} \frac{h_m}{(2m+1)^n} + ... + 2 \sum^{k}_{m=1} \frac{h^{(n)}_m}{(2m+1)}~,
\end{eqnarray}
or
\begin{eqnarray}
F(k) = \sum^{k}_{m=1} \frac{H_{m-1}}{(2m-1)^n} +2\sum^{k}_{m=1} \frac{h_m}{(2m-1)^n} + ... + 2\sum^{k}_{m=1} \frac{h^{(n)}_m}{(2m-1)}
- 2n h^{(n+1)}_k~.
\end{eqnarray}
Therefore we get:
\begin{eqnarray}
\sum^{k-1}_{m=1}\frac{h^{(n)}_{m}}{m(m-k)} = \frac{1}{k} \sum^{k-1}_{m=1}\frac{h^{(n)}_m}{m} +\frac{1}{k} \sum^{k}_{m=1} \frac{H_{m-1}}{(2m-1)^n}
+\frac{2}{k}\sum^{k}_{m=1} \frac{h_m}{(2m-1)^n} + ... + \frac{2}{k}\sum^{k}_{m=1} \frac{h^{(n)}_m}{(2m-1)} - \frac{2n}{k} h^{(n+1)}_k~.
\end{eqnarray}
Finally Eq.~(90) results in: 
\begin{eqnarray}
f_n(k)&=& \frac{1}{k} \sum^{k-1}_{m=1}\frac{h^{(n)}_m}{m} +\frac{1}{k} \sum^{k}_{m=1} \frac{H_{m-1}}{(2m-1)^n} +
\frac{2}{k}\sum^{k}_{m=1} \frac{h_m}{(2m-1)^n} + ... + \frac{2}{k}\sum^{k}_{m=1} \frac{h^{(n)}_m}{(2m-1)} - \frac{2n}{k} h^{(n+1)}_k
\nonumber \\ &-& \sum^{k}_{i=1}\frac{h^{(n)}_i}{i} - \frac{1}{k} \sum^{k}_{m=1} \left( \sum^{i}_{m=1} \frac{1}{i(2i-2m+1)^n} \right)~.
\end{eqnarray}
By partial fraction decomposition of the last inner sum it follows:
\begin{eqnarray}
f_n(k)=\frac{1}{k}\sum^{k}_{m=1}\frac{H_{m-1}}{(2m-1)^n}-\frac{h^{(n)}_k}{k^2}+2ln(2)\frac{h^{(n)}_k}{k}+\frac{3}{2}\zeta(2) \frac{h^{(n-1)}_k}{k}
+ ... +2\left(1-\frac{1}{2^{n}} \right) \zeta(2n)\frac{h^{(1)}_k}{k}-2n \frac{h^{(n+1)}_k}{k}~.
\end{eqnarray}
It remains to calculate 
\begin{eqnarray}
\sum^{\infty}_{i=1}\frac{h^{(n)}_{i}}{i(k+i)} = \frac{1}{k} \sum^{\infty}_{i=1} \frac{1}{i(2i-1)^{2n}} + \frac{1}{k}
\sum^{k-1}_{m=1} \left( \sum^{\infty}_{i=1} \frac{1}{(i+m)(2i-1)^n} \right)~.
\end{eqnarray}
By partial fraction decomposition of the last sum on right side it follows:
\begin{eqnarray}
\sum^{\infty}_{i=1}\frac{h^{(n)}_{i}}{i(k+i)} &=& \frac{1}{k} \sum^{\infty}_{i=1} \frac{1}{i(2i-1)^{2n}} + \frac{3}{2}\zeta(2)\frac{1}{k}
\left( h^{(n-1)}_i - 1 \right) + ... +2(-)^n \left(1-\frac{1}{2^{n}} \right)\zeta(2n) \left( h^{(1)}_i - 1 \right) \nonumber \\ &-&
\frac{1}{k} \sum^{k-1}_{i=1} \frac{H_{i}}{(2i+1)^{2n}} - 2ln(2)\frac{1}{k}\sum^{k-1}_{i=1}\frac{1}{(2i+1)^n}~,
\end{eqnarray}
and finally we get:
\begin{eqnarray}
\sum^{\infty}_{i=1}\frac{h^{(n)}_{i}}{i(k+i)} = \frac{3}{2}\zeta(2)\frac{1}{k} h^{(n-1)}_i + ... +2(-)^n \frac{1}{2^{n}} \zeta(2n) h^{(1)}_i
-\frac{1}{k} \sum^{k}_{i=1} \frac{H_{i-1}}{(2i-1)^{2n}} - 2ln(2)\frac{h^{(n)}_k}{k}~.
\end{eqnarray}
Inserting Eq.~(103) in Eq.~(101) and separating in odd and even contributions Eq.~(88) and Eq.~(89) follow. Thus lemma 3 is proofed.

As a first example we compute Eq.~(89) for n=1.
\begin{eqnarray}
\sum^{\infty}_{k=1}\frac{1}{k^3}\sum^{\infty}_{i=1,i\ne k}\frac{h^{(2)}_i}{i(k-i)} = -\sum^{\infty}_{k=1}\frac{1}{k^3}
\sum^{\infty}_{i=1}\frac{h^{(2)}_i}{i(k+i)} + 3\zeta(2)\sum^{\infty}_{k=1}\frac{h_k}{k^4}-
\sum^{\infty}_{k=1}\frac{h_k^{(2)}}{k^5} -4\sum^{\infty}_{k=1}\frac{h_k^{(3)}}{k^4} .
\end{eqnarray}
First we get:
\begin{eqnarray}
\sum^{\infty}_{i=1}\frac{h^{(2)}_i}{i}\sum^{\infty}_{k=1,k\ne i}\frac{1}{k^3(i-k)} = 
\sum^{\infty}_{i=1}\frac{h^{(2)}_i}{i} \left(\frac{H_i}{i^3}-\frac{4}{i^4}+\zeta(3)\frac{1}{i}+\zeta(2)\frac{1}{i^2} \right) \nonumber
\end{eqnarray}
\begin{eqnarray}
= \sum^{\infty}_{i=1}\frac{H_i h^{(2)}_i}{i^4} - 4 \sum^{\infty}_{i=1}\frac{h^{(2)}_i}{i^5} + \zeta(3)
\sum^{\infty}_{i=1}\frac{h^{(2)}_i}{i^2} + \zeta(2) \sum^{\infty}_{i=1}\frac{h^{(2)}_i}{i^3}~.
\end{eqnarray}
Furthermore it follows:
\begin{eqnarray}
\sum^{\infty}_{k=1}\frac{1}{k^3} \sum^{\infty}_{i=1}\frac{h^{(2)}_i}{i(k+i)} = 
\sum^{\infty}_{i=1}\frac{h^{(2)}_i}{i}\sum^{\infty}_{k=1}\frac{1}{k^3(i+k)} = \zeta(3)\sum^{\infty}_{i=1}\frac{h^{(2)}_i}{i^2}
-\zeta(2)\sum^{\infty}_{i=1}\frac{h^{(2)}_i}{i^3} + \sum^{\infty}_{i=1}\frac{H_i h^{(2)}_i}{i^4}~.
\end{eqnarray}
Inserting this expression in Eq.~(101) and comparing with Eq.~(104) we get:
\begin{eqnarray}
\sum^{\infty}_{k=1}\frac{h_k^{(3)}}{k^4} = \frac{3}{4}\zeta(2)\sum^{\infty}_{k=1}\frac{h_k}{k^4} + 
\frac{1}{2}\zeta(2)\sum^{\infty}_{k=1}\frac{h_k^{(2)}}{k^3} -\frac{5}{4}\sum^{\infty}_{k=1}\frac{h_k^{(2)}}{k^5}~.
\end{eqnarray}
All sums on the right side of Eq.~(107) are known. This way Eq.~(86) results. The other five subtypes $T^2_{3,m}, ..., T^6_{3,m}$ 
belonging to this family of Euler sums follow analogously either using the corresponding identities formulated by Adegoke or from
the corresponding calculational procedures introduced in this work.

A second example may be given for n=3, where we compute $T^1_{4,5}$. It follows first:
\begin{eqnarray}
\sum^{\infty}_{k=1}\frac{1}{k^4}\sum^{\infty}_{i=1,i\ne k}\frac{h^{(3)}_i}{i(k-i)} = \sum^{\infty}_{k=1}\frac{1}{k^4}
\sum^{\infty}_{i=1}\frac{h^{(3)}_i}{i(k+i)} + 3\zeta(2)\sum^{\infty}_{k=1}\frac{h_k^{(2)}}{k^5}-
\sum^{\infty}_{k=1}\frac{h_k^{(3)}}{k^6} -6 \sum^{\infty}_{k=1}\frac{h_k^{(4)}}{k^5}~.
\end{eqnarray}
With
\begin{eqnarray}
\sum^{\infty}_{i=1}\frac{h^{(3)}_i}{i}\sum^{\infty}_{k=1,k\ne i}\frac{1}{k^4(i-k)} = 
-\sum^{\infty}_{i=1}\frac{h^{(3)}_i}{i} \left(\frac{H_i}{i^4}-\frac{5}{i^5} +\zeta(4)\frac{1}{i} +\zeta(3)\frac{1}{i^2}
+\zeta(2)\frac{1}{i^3} \right) \nonumber
\end{eqnarray}
\begin{eqnarray}
= -\sum^{\infty}_{i=1}\frac{H_i h^{(3)}_i}{i^5} + 5 \sum^{\infty}_{i=1}\frac{h^{(3)}_i}{i^6}
- \zeta(4) \sum^{\infty}_{i=1}\frac{h^{(3)}_i}{i^2} -\zeta(3) \sum^{\infty}_{i=1}\frac{h^{(3)}_i}{i^3} + 
\zeta(2) \sum^{\infty}_{i=1}\frac{h^{(3)}_i}{i^4}~. 
\end{eqnarray}
Furthermore it follows:
\begin{eqnarray}
\sum^{\infty}_{k=1}\frac{1}{k^4} \sum^{\infty}_{i=1}\frac{h^{(3)}_i}{i(k+i)} = 
\sum^{\infty}_{i=1}\frac{h^{(3)}_i}{i}\sum^{\infty}_{k=1}\frac{1}{k^4(i+k)} = \zeta(4)\sum^{\infty}_{i=1}\frac{h^{(3)}_i}{i^2} 
- \zeta(3)\sum^{\infty}_{i=1}\frac{h^{(3)}_i}{i^3}
+\zeta(2)\sum^{\infty}_{i=1}\frac{h^{(3)}_i}{i^4} - \sum^{\infty}_{i=1}\frac{H_i h^{(3)}_i}{i^5}~.
\end{eqnarray}
Inserting this expression in Eq.~(108) and comparing with Eq.~(109) we get:
\begin{eqnarray}
\sum^{\infty}_{k=1}\frac{h_k^{(4)}}{k^5} = \frac{1}{3}\zeta(4)\sum^{\infty}_{k=1}\frac{h_k^{(3)}}{k^2} + 
\frac{1}{2}\zeta(2)\sum^{\infty}_{k=1}\frac{h_k^{(2)}}{k^5} + \frac{1}{3}\zeta(2)\sum^{\infty}_{k=1}\frac{h_k^{(3)}}{k^4} 
-\sum^{\infty}_{k=1}\frac{h_k^{(3)}}{k^6} .
\end{eqnarray}
Again all sums on the right side of Eq.~(111) are known. Thus $T^1_{4,5}$ follows explicitly in terms of zeta values only. 
It follows numerically $T^1_{4,5} = 1.0373935033868233$ with an accuracy of $10^{-16}$ compared to the analytical expression.

\section{Summary}
We have introduced a special summation method that allows, based on proper two-valued integer functions, to calculate explicitly
in terms of zeta values a variety of Euler sums involving odd-type harmonic numbers.


\begin{thebibliography}{99}
\bibitem{fla97}
Philippe Flajolet and Bruno Salvey, Euler sums and Contour integral representations, Experimental Mathematics 7 (1997) 1, 15.

\bibitem{sit85}
R. Sitaramachandrarao, A Formula of S. Ramanujan, Journal of number theory 25, (1987) 1-19.

\bibitem{zeh07}
De-Yin Zheng, Further summation formulae related to generalized harmonic numbers, J. Math. Anal. Appl. 335 (2007) 692–70.

\bibitem{ade16}
Kunle Adegoke1 and Olawanle Layeni, New Finite and Infinite Summation Identities Involving the Generalized Harmonic Numbers,
J. Ana. Num. Theor.4, No. 1, (2016) 49-60.

\bibitem{ce17}
Ce Xu, Yingyue Yang and Jianwen Zhang, Explicit evaluation of quadratic Euler sums,  International Journal of Number Theory, Vol. 13,
No. 03, (2017) 655-672. 

\bibitem{ce18}
Ce Xu, EULER SUMS OF GENERALIZED HYPERHARMONIC NUMBERS, Korean Math. Soc.55, No. 5, (2018) 1207–1220.
 
\bibitem{ce20}
Ce Xu, Evaluations of Euler-Type Sums of Weight $\le$ 5, Bulletin of the Malaysian Mathematical Sciences Society volume 43, (2020)
847–877. 

\bibitem{bra20}
J. Braun, D. Romberger and H. J. Bentz, On four families of power series involving harmonic numbersand central binomial coefficients,
arxiv:2006.13115v1 (2020).

\end{thebibliography}
\end{document}